\documentclass[12pt,a4paper]{amsart}
\usepackage{amssymb,amsmath,amsthm}
\usepackage{graphicx}
\usepackage{enumerate}
\usepackage{tikz}

\long\def\onefigure#1#2{
\begin{figure*}[tbp]
\begin{center}
#1
\end{center}
\caption{#2}
\end{figure*}
} 

\newcommand{\lipefig}[2]  
{\onefigure{\mbox{\psfig{file=#1.eps}}}{\label{f:#1} #2} }

\newtheorem{theorem}{Theorem}[section]
\newtheorem{lemma}[theorem]{Lemma}

\newtheorem{fact}{Fact}[section]

\newcommand{\remove}[1]{}

\newcommand{\de}{\delta}

\newcommand{\la}{\lambda}

\newcommand{\eps}{\varepsilon}

\newcommand{\rr}{\mathbb{R}^d}
\newcommand{\R}{\mathbb{R}}
\newcommand{\gr}{\mathrm{Gr}}

\newcommand{\F}{\mathcal{F}}

\newcommand{\G}{\mathcal{G}}
\newcommand{\inn}{\mathrm{int\;}}
\newcommand{\conv}{\mathrm{conv}}
\newcommand{\pos}{\mathrm{pos\;}}

\newcommand{\bd}{\mathrm{bd\;}}

\newcommand{\aff}{\textrm{aff\;}}
\newcommand{\lin}{\textrm{lin\;}}

\newcommand{\relint}{\textrm{relint}}

\numberwithin{equation}{section}

\begin{document}

\title{Erd\H os-Szekeres theorem for $k$-flats}
\author{Imre B\'ar\'any, Gil Kalai, and Attila P\'or}


\begin{abstract} We extend the famous Erd\H os-Szekeres theorem to $k$-flats in $\rr$.
\end{abstract}

\maketitle

\section{Introduction}\label{sec:introd}

The famous Erd\H os-Szekeres \cite{ESz35} theorem from 1935 asserts that for every $n\ge 3$ there is an integer $N=N(n)$ such that any set of $N$ points in general position in the plane contains a subset of size $n$ in convex position, that is, these points form the vertices of a convex $n$-gon. General position means that no three points are collinear. It is trivial yet important to note that not every set of $n>3$ points is in convex position. In addition, the vertices of a convex $N$-gon show that being convex is the only type of $n$-tuples that appears in every set of size $N$ in general position.

\smallskip
We want to extend the Erd\H os-Szekeres theorem to $k$-flats (affine subspaces of dimension $k$) in $\rr$. We say that a set of $n$ $k$-flats in $\rr$, $\{U_1,\ldots,U_n\}$, is in {\sl convex position}, or simply that it is {\sl convex} if there is a $d$-dimensional polytope $P \subset \rr$ such that $U_i \cap P$ is a $k$-dimensional face of $P$ for every $i\in [n]:=\{1,\ldots,n\}$. Here of course $0\le k\le d-1$. This definition is a direct extension of the Erd\H os-Szekeres setting which corresponds to the case $k=0$.

\smallskip
We need a general position assumption for $n$-tuples of $k$-flats. Suppose $d\ge 3$ and $0\le k \le d-2$. An $n$-tuple $U_1,\ldots,U_n$ of $k$-flats is in {\sl general position} if there is a $(d-k)$-flat $A$ such that $U_i\cap A$ is a single point for every $i\in [n]$, no three of these $n$ points are collinear, and their affine hull coincides with $A$. We note that for $k=0$ this gives a weaker condition than
what is typically considered general position for points in $\rr$.
The case of hyperplanes, that is $k=d-1$, is different: $n\ge d$ hyperplanes in $\rr$ are in general position if every $d$ of them intersect in a single point and these ${n \choose d}$ points are distinct.

\begin{theorem}\label{th:k-flats} Assume $d\ge 2,k,n$ are integers with $0\le k \le d-1$ and $n\ge d-k+1$. There is an integer $N=N(k,d,n)$ such that every set of $N$ $k$-flats in $\rr$ in general position contains a convex subset of size $n$.
\end{theorem}

This is the first main result of our paper. The second one says that, for large enough $n$, there are non-convex sets of $k$-flats of size $n$.

\begin{theorem}\label{th:non-conv} Assume $d\ge 2,k$ are integers with $0\le k \le d-1$. Then there is an integer $n=n(d,k)$ and there is an $n$-tuple of $k$-flats in $\rr$ in general position which is not convex.
\end{theorem}

\smallskip
We mention that the same definition of convex position for $k$-flats is given in \cite{BBLP}. The target there is different. Namely it is shown in \cite{BBLP} that, given a finite family of lines in $\R^2$ in a suitably general position, if every $5$-tuple of lines in the family is in convex position, then all of them are in convex position. It is simple to see that every $4$-tuple of lines is in convex position. So the result in \cite{BBLP} says that this type of Helly number of lines in $\R^2$ is $5$. It is not clear if the analogous Helly number for $k$-flats in $\rr$ is finite or not, and if it is, what its value is. 

\smallskip
Rich collections of results and questions around the Erd\H os-Szekeres theorem can be found in \cite{MorSol}, \cite{BarKar}, and \cite{MorSol2}.

\bigskip
\section{Background and preparations}\label{sec:backg}

Define $N=N(n)$ as the smallest integer for which the Erd\H os-Szekeres theorem holds. The bound $N(n) \le {2n-4 \choose n-2}+1$ is from their 1935 paper~\cite{ESz35}. 
A recent breakthrough by Suk~\cite{Su} asserting that $N(n)\le 2^{n+o(n)}$ comes close to the so called happy-end conjecture saying that $N(n)=2^{n-2}+1$. There is actually little evidence supporting this conjecture except the lower bound $N(n) \ge 2^{n-2}+1$ from~\cite{ESz61}.

\smallskip
There is an alternative way for defining convex position of $k$-flats:

\begin{lemma}\label{l:equiv} Assume $\{U_1,\ldots,U_n\}$ is a set of $k$-flats in $\rr$. They are in convex position if and only if there is a $d$-dimensional strictly convex body $K$ such that $K\cap U_i$ is a single point for every $i \in [n]$.
\end{lemma}

{\bf Proof.} Assume first that these $k$-flats are in convex position, that is, there is a  $d$-dimensional polytope $P \subset \rr$ such that $P\cap U_i$ is a  $k$-dimensional face of $P$ for every $i \in [n]$. Let $h_i$ be a hyperplane with $h_i\cap P=U_i\cap P$ and let $a_i$ be the center of gravity of the face  $U_i\cap P$. Fix a large ball $B_i$ such that $h_i$ is tangent to $B_i$ at $a_i$ and with radius so large that each $a_j,\;(j\ne i)$ is contained in the interior of $B_i$. Then $\bigcap_1^n B_i$ is a ($d$-dimensional) convex body in $\rr$ with the required properties.

Conversely, let $K$ be a strictly convex body with $U_i\cap K$ a single point, $a_i$ say, for all $i\in [n]$. Let $h_i$ be a hyperplane containing $U_i$ and intersecting $K$ in only the point $a_i$. There is a small ball $B_i^*$ centered at $a_i$ such that for every distinct $i,j \in [n]$, ($i \neq j$), 
the ball $B_i^*$ and the hyperplane $h_j$ are disjoint. 
Now choose a $k$-dimensional polytope $P_i$ lying in $U_i \cap B_i^*$ and $a_i \in P_i$. Finally let $P_0$ be a $d$-dimensional polytope contained in $\inn K$. Then $P=\conv \bigcup_0^n P_i$ is a polytope with the required properties. $P_0$ was needed to make sure that $P$ is $d$-dimensional.  \qed

\smallskip
We fix some notation. Define $\conv X$, $\pos X$, $\lin X$, and $\aff X$ as the convex, cone, linear, and affine hull of a set $X\subset \rr$. We write $uv$ for the scalar product of vectors $u,v \in \rr$.
Let $B$ denote the Euclidean unit ball of $\rr$ centered at the origin.

\smallskip
We are going to work with $\gr(k,d)$, the Grassmannian of $k$-dimensional subspaces of $\rr$. The Grassmannian is a metric space where the distance $d(U,V)$ of $U,V \in \gr(k,d)$ is given via
\[
d(U,V) \mbox{ is the Haussdorff distance of the sets } U\cap B \mbox{ and }V\cap B.
\]

If $U$ and $V$ are $k$-dimensional subspaces in $\rr$,
the angle between $U$ and $V$ is defined 
by C. Jordan \cite{Jor1875} as
as the largest of the principal angles 
\[
\angle(U,V)= \arcsin (d(U,V)),
\]
We have the following facts
\[
\frac 2{\pi}\angle(U,V) \le d(U,V) \le \angle(U,V),
\]
and for every unit vector $u \in U$ there exists a unit vector $v \in V$ such that
$\angle (u,v) \le \angle (U,V)$ which implies
$\|u-v\| \le \angle(U,V)$.

Then $\gr(k,d)$ is a compact metric space so for every $\eps>0$ it contains a finite set $\G_{\eps}$, called an $\eps$-net, whose size depends only on $\eps,k,d$ such that for every $U \in \gr(k,d)$ there is a $V \in \G_{\eps}$ such that $\angle(U,V)<\eps$.

We are going prove the following result which is the cone-version of Theorem~\ref{th:non-conv}.

\begin{theorem}\label{th:non-cone} Assume $k,d$ are integers with $0<k<d$. 
	Then for some $\eps=\eps(k,d)>0$ any $\eps$-net $\G_{\eps}$ has the property that there is no polyhedral cone $C \subset \rr$ such that $U\cap C$ is a $k$-dimensional face of $C$ for every $U \in \G_{\eps}$.
\end{theorem}

\bigskip
\section{Proof of Theorem~\ref{th:k-flats}}\label{sec:proofofmain}

The following simple lemma shows that the Erd\H os-Szekeres theorem implies the case $k=0$ of Theorem~\ref{th:k-flats}.

\begin{lemma}\label{l:casek=0} $N(0,d,n)\le N(n)$.
\end{lemma}

{\bf Proof.} This is a well known argument. Project the $N=N(n)$ points in $\rr$ to a 2-dimensional plane $L$ which can be chosen so that no three of the projected points are collinear because of the general position assumption on the $N$ points in $\R^d$. The Erd\H os-Szekeres theorem implies then that there is an $n$-set in convex position among the $N$ projected points. It is clear that the corresponding set of $n$ points in $\rr$ is in convex position as well. \qed

\smallskip
Suk's result cited above implies that $N(0,d,n)\le N(n) \le 2^{n+o(n)}$. We mention a brave conjecture of Zolt\'an F\"uredi (under a slightly stronger general position assumption, see \cite{KarVal}): $N(0,d,n)=O(2^{n^{1/(d-1)}})$ which looks very difficult and which is supported by a lower bound of the same order~\cite{KarVal}.

\smallskip
The case $d=2,k=1$ of Theorem~\ref{th:k-flats} is about $N(1,2,n)$ lines in the plane. The main theorem in \cite{BRT} says that among $N={2n-4 \choose n-2}$ lines in the plane in general position there are always $n$ convex. Here general position means that no two lines are parallel and no three lines are concurrent which is the same in the planar case as our general position condition for hyperplanes in $\rr$. So $N(1,2,n) \le {2n-4 \choose n-2}$ which is of order $4^n/\sqrt n$. The paper~\cite{BRT} also gives a lower bound of order $4^n/n$.  The boundedness of $N(d-1,d,n)$ follows from the planar case as the following lemma shows.

\begin{lemma}\label{l:casek=d-1} $N(d-1,d,n)\le N(1,2,n)$.
\end{lemma}

{\bf Proof.} Given $N=N(1,2,n)$ hyperplanes $H_1,\ldots,H_N$ in $\rr$ in general position, there is a 2-dimensional plane $L$ such that the intersections $H_i\cap L$ are lines in $L$ in general position. This can be checked easily, details are left to the interested reader. Among these $N$ lines in $L$ there are $n$ in convex position. The complement of the corresponding $n$ hyperplanes, $H_{i_1},\ldots,H_{i_n}$, in $\rr$ consists of finitely many connected components. Each such component is the intersection of finitely many (at most $n$) open halfspaces, so is an open polyhedron. The intersection of one of them, $C$ say, with $L$ is a (possibly unbounded) open convex $n$-gon. It is clear that the closure of $C$ is a $d$-dimensional polyhedron $Q$ that has $n$ facets, each one of the form $H_{i_j}\cap Q.$ It is evident that $Q$ contains a polytope $P$ such that $H_{i_j}\cap P$ is a facet of $P$ for every $j \in [n].$ So these $n$ hyperplanes are in convex position. \qed

\smallskip
The proof shows the upper bound $N(d-1,d,n)\le {2n-4 \choose n-2}$. It is not clear how good this upper bound is.

\medskip
{\bf Proof} of Theorem~\ref{th:k-flats}. In view of the previous lemmas we can assume that $0<k<d-1$ and so $d\ge 3$. Set $N=N(0,d-k,n)$ which is finite because of Lemma~\ref{l:casek=0}, and consider a set of $k$-flats $U_1,\ldots,U_N$ in $\rr$ in general position. Then there is a $(d-k)$-flat $A$ such that $a_i:=A \cap U_i$ is a single point for all $i\in [N]$, no three points from the set $X=\{a_1,\ldots,a_N\}$ are collinear.
By Lemma~\ref{l:casek=0} again $X$ contains an $n$-element subset $Y=\{b_1,\ldots,b_n\}$ that forms the vertex set of a convex polytope $Q$.
For simpler notation we assume that $b_i = a_i$.
The dimension of $A^* = \aff Y$ is at least $2$ and at most $d-k$.
Let $h_i^*$ be a hyperplane in $A^*$ tangent to $Q$ at $b_i$, that is $h_i^* \cap Q = \{ b_i\}$.
The hyperplane $h_i = \aff ( h_i^* \cup U_i)$ in $\rr$ is disjoint from
$\conv (Y \setminus \{ b_i\} )$. 
Let $b_0 \in \relint Q$ be a point in the relative interior of $Q$.
Then for some small $\de>0$, for every $0\le j \le n$, $1\le i \le n$ and $i \neq j$
the ball $b_j+\de B$ is disjoint from $h_i$.
For every $1\le j \le n$ choose a $k$-dimensional polytope $Q_j$ in $(b_j+\de B) \cap U_j$
and choose $Q_0$ to be a $d$-dimensional polytope in  $b_0+\de B$.
Set $P=\conv (\bigcup_{0\le j \le n} Q_j)$. Then $P$ is a $d$-dimensional polytope in $\rr$ and for every $1 \le j \le n$ each $Q_j = P \cap U_j$ is a $k$-face of $P$. \qed

\smallskip
This implies that in the range $0<k<d-1$ (and $d\ge 3$), $N(k,d,n)\le N(0,d-k,n)\le 2^{n+o(n)} $. F\"uredi's conjecture, if true, would imply much better upper bounds. We wonder for instance what the value of $N(d/2,d,n)$ could be. 

\bigskip
\section{Proof of Theorem~\ref{th:non-cone}}\label{sec:non-cone}

We begin by assuming that for some $\eps>0$ and for some $\eps$-net $\G_{\eps}\subset \gr(k,d)$ there is a polyhedral cone $C$ in $\rr$ such that $U\cap C$ is a $k$-face of $C$ for every $U\in \G_{\eps}$, and show that $\eps$ has to be larger than some positive constant that only depends on $k$ and $d$.

\smallskip
Under the above assumption for every $U \in \G_{\eps}$ we must have
\begin{enumerate}[{(1)}]
 \item  $U\cap C\ne \{0\}$,  and
 \item  $U \cap \inn C = \emptyset$.
\end{enumerate}
where $\inn C$ denotes the interior of $C$.

\smallskip
We need the following facts.

\begin{fact}\label{f:angles} For every $V \in\gr(s,d)$ with $s\le d-k$ there is $U \in \G_{\eps}$ such that $|uv|< \eps$ for every pair of unit vectors $u \in U$ and $v\in V$.
\end{fact}

The {\bf proof} is simple. The orthogonal complement of $V\in \gr(s,d)$ is of dimension $d-s\ge k$ so it contains a subspace $U_0$ from $\gr(k,d)$. Then $\G_{\eps}$ contains a subspace $U$ with $\angle (U,U_0)< \eps$. For a unit vector $u \in U$ there is a unit vector $u_0\in U_0$ with $\|u-u_0\|<\eps$. Then for a unit vector $v \in V$ $|uv|=|(u-u_0)v+u_0v|=|(u-u_0)v|<\eps$. \qed

\begin{fact}\label{f:determinant} If $M$ is a $t \times t$ matrix with all diagonal entries $1$ and every other entry at most $\de$ in absolute value, then
$\det M \ge 1-t!\de$.
\end{fact}

The {\bf proof} follows from the Leibniz formula as the product of the diagonal entries 
is $1$, and every other term is at most $\de$ in absolute value.

\smallskip
We are going to construct unit vectors $a_0,a_1,\ldots,a_{d-k}$ on the boundary, $\bd C$, of $C$ that are pairwise {\sl almost orthogonal} meaning that $|a_ia_j|<\eps$ for distinct $i,j$.
We begin by selecting a unit vector $a_0\in \bd C$. Then $V_1=\lin \{a_0\}$ is an element in $\gr(1,d)$ so by Fact~\ref{f:angles} there is $U_1\in \G_{\eps}$ with $|ua_0|<\eps$ for every unit vector $u \in U_1$. In view of conditions (1) and (2) we can choose a unit vector $a_1\in U_1\cap \bd C$. Consequently 
$|a_0a_1|< \eps$. Assume that for $j \le d-k$ we have unit vectors $a_0,\ldots,a_{j-1} \in \bd C$ such that 
$|a_ia_h|<\eps$ for distinct $i,h \in \{0,1,\ldots,j-1\}$. The subspace $V_j:=\lin \{a_0,\ldots,a_{j-1}\} \in \gr(j,d)$ is of dimension $j \le d-k$, and again by Fact~\ref{f:angles}, there is $U_j \in \G_{\eps}$ such that $|uv|<\eps$ for every pair of unit vectors $u\in U_j$ and $v\in V_j$. Choosing a unit vector $a_j \in U_j \cap \bd C$, again by conditions (1) and (2), finishes the construction.

\smallskip
Next we find linearly independent unit vectors $b_i \in \inn C$ very close to $a_i$ (for every $i\in \{0,1,\ldots,d-k\}$) so that $|b_ib_h|<\eps$ for distinct $i,h$. This is clearly possible. 
The $(d-k+1)$-dimensional cone   $D:=\pos \{b_0,b_1,\ldots,b_{d-k}\}$ lies in the subspace $V= \lin \{b_0,b_1,\ldots,b_{d-k}\}\in \gr(d-k+1,d)$. Our target is to show that for some $U \in \G_{\eps}$ the intersection $D\cap U$ is a halfline which would contradict condition (2) because $D$ lies in the interior of $C$.

\smallskip
Let $c_1,\ldots,c_{k-1}$ be an orthonormal basis of 
the orthogonal complement of $V$ and define $W=\lin \{c_1,\ldots,c_{k-1},b\}\in \gr(k,d)$ where $b=b_0+b_1+\ldots +b_{d-k}$. Consider the linear system of equations
\begin{equation}\label{eq:sys}
\sum_0^{d-k}x_ib_i+\sum_1^{k-1}y_jc_j-yb=0.
\end{equation}

\begin{lemma}\label{f:syst}  If $\eps>0$ is small enough, then the only solution to the system (\ref{eq:sys}) is $x_0=\ldots =x_{d-k}=1$, $y_1=\ldots =y_{k-1}=0, y=1$ and its scalar multiples.
\end{lemma}

{\bf Proof.} Let $M$ be the matrix with columns $b_0,\ldots,b_{d-k},c_1,\ldots,c_{k-1}$. So $M$ is a $d\times d$ matrix. As the system (4.1) is homogeneous, it is enough to check that $\det M\ne 0$ or, what is the same, $\det M^TM \ne 0$. Every entry on the main diagonal of $M^TM$ is 1, all other entries are at most $\eps$ in absolute value.  
By Fact~\ref{f:determinant}  $\det M^TM>1-d!\eps>\frac 14$ if $\eps< \frac 3{4d!}$. \qed

\smallskip
We assume now that $\det M>\frac 12$ (replace $c_1$ by $-c_1$ if $\det M<0$). Let $U \in \G_{\eps}$ be a subspace with $\angle (U,W)< \eps$ and choose vectors $c_1^*,\ldots,c_{k-1}^*,b^*$ of $U$ with $\|c_j-c_j^*\|< \eps$ and $\|c_j^*\|=1$ for all $j\in [k-1]$ and $\angle (b,b^*)<\eps$ and $\|b\|=\|b^*\|$. Consider the system
\begin{equation}\label{eq:sys*}
\sum_0^{d-k}x_ib_i+\sum_1^{k-1}y_jc_j^*-yb^*=0.
\end{equation}

\begin{lemma}\label{l:syst*}  If $\eps>0$ is small enough, then the system {\rm (4.2)} has a solution with $x_i>0$ for all $i\in \{0,1,\ldots,d-k\}$.
\end{lemma}

{\bf Proof.} Let $M^*$ be the matrix with columns $b_0,\ldots,b_{d-k},c_1^*,\ldots,c_{k-1}^*$. We check first $\det M^{*T}M^{*}>\frac 14$ if $\eps<\frac 1{4d!}$. All entries on the main diagonal are equal to one, the entries $b_ib_h$ are at most $\eps$ in absolute value. The entry $b_ic_j^*=b_i(c_j^*-c_j)+b_ic_j=b_i(c_j^*-c_j)$ so $|b_ic_j^*| <\eps$. Finally for distinct $j,h$ 
\[ 
c_j^*c_h^*=(c_j^*-c_j)(c_h^*-c_h)+c_j^*c_h+c_jc_h^*-c_jc_h
\] 
so $|c_j^*c_h^*|<\eps^2+\eps+\eps+0<3\eps$. 
Again by Fact~\ref{f:determinant}
$\det M^{*T}M^{*}>1-d!3\eps>\frac 14$ if $\eps< \frac 1{4d!}$. This implies that the solution to (\ref{eq:sys*}) is unique up to a scalar multiplier. Fix now such a solution so that the maximal absolute value of the numbers in the set $S=\{x_0,\ldots,x_{d-k},y_1,\ldots,y_{k-1},y\}$ equals one. 

\smallskip
Observe that $c_j^*c_i=(c_j^*-c_j)c_i+c_jc_i$ and $c_jc_i=0$ if $i\ne j$ and $c_jc_i=1$ if $i= j$. Similarly $b^*c_i=(b^*-b)c_i+bc_i$ and here $bc_i=0$, so $|b^*c_i|\le \|b^*-b\|\le \eps \|b\|\le \eps (d-k+1)$. Multiplying equation (\ref{eq:sys*}) by $c_i$ gives $y_i+\sum_1^{k-1} y_j(c_j^*-c_j)c_i-y(b^*-b)c_i=0$ implying
\begin{eqnarray*}
|y_i|&=& \left |\sum_1^{k-1} y_j(c_j^*-c_j)c_i-y(b^*-b)c_i \right |\\
      &\le& \sum_1^{k-1} |y_j|\|c_j^*-c_j\|+\eps |y|(d-k+1)<\eps d,
\end{eqnarray*}
because $|y_j|\le 1$ and $|y|\le 1$. Then each $|y_i|< \frac 14$ if $\eps \le \frac 1{4d}$ and none of the $y_j$ is maximal in $S$.

\smallskip
Multiply equation (\ref{eq:sys*}) by $b_i$. Again $c_j^*b_i=(c_j^*-c_j)b_i+c_jb_i=(c_j^*-c_j)b_i$ showing that $|c_j^*b_i|<\eps$. Analogously $b^*b_i=(b^*-b)b_i+bb_i=(b^*-b)b_i+1+\sum_{j\neq i} b_ib_j$ 
and here $|(b^*-b)b_i|\le \eps\|b\|\le \eps(d-k+1)$. Thus we have $\sum x_j b_ib_j+\sum_1^{k-1}y_j(c_j^*-c_j)b_i- y\left[(b^*-b)b_i+1+\sum_{j\neq i} b_ib_j\right]=0$ showing that
\begin{eqnarray*}
|x_i-y|&=&\left |\sum_1^{k-1}y_j(c_j^*-c_j)b_i+\sum_{j\neq i} x_jb_ib_j
- y\left [(b^*-b)b_i +\sum_{j\neq i} b_ib_j \right ] \right |\\
          &<& \eps \sum_1^{k-1}|y_j|+ \eps (d-k)
          +\eps |y|(d-k+1) + \eps|y| (d-k) \\
          &\le& \eps (3d-2k) < 3\eps d.
\end{eqnarray*} 
Thus $x_i$ and $y$ differ by at most $\frac 14$ if $\eps <\frac 1{12d}$. Now either some $x_i$ or $y$ has maximal absolute value in $S$ equal to one. We can assume that either $x_i=1$ for some $i$ or $y=1$ (by multiplying the solution by $-1$ if necessary). In either case  $x_j>\frac 12$ for all $j=0,1,\ldots,d-k+1$.

\smallskip
Then $z=\sum_0^{d-k}x_ib_i\in \inn D$ and $z=yb^*-\sum_1^{k-1}y_jc_j^*\in U$. Consequently $z$ is a common point of $\inn D$ and of $U \in \G_{\eps}$ provided $\eps< \min \{ \frac 1{12d}, \frac 1{4d!} \}$ contradicting condition (2).
\qed

\section{Proof of Theorem~\ref{th:non-conv}}\label{sec:non-conv}

This follows directly from the example for the cone version. Indeed, consider the example given in Theorem~\ref{th:non-cone} of an $\eps$-net $\G_{\eps}$ in $\gr(k+1,d+1)$ and let $H$ be a hyperplane in general position with respect to $\G_{\eps}$ with $0 \notin H$. General position simply means that $U\cap H $ is a $k$-flat in $H$ for every $U \in \G_{\eps}$. Then $H$ can be taken for $\rr$ and the system of $k$-flats $\{U \cap H:U \in \G_{\eps}\}$ is not in convex position. Because if it were and $P$ were the polytope in $H$ such that $U\cap H \cap P$ is a $k$-face of $P$ for every $U \in \G_{\eps}$, then the cone $\pos P$ in $\R^{d+1}$ would show that the subspaces in $\G_{\eps}$ are in convex position. \qed

\smallskip
In conclusion we give a simple example of a non-convex family $\F$ of $2^d+d$
hyperplanes in $\rr$. First let $H_i$ be the hyperplane with equation $x_i=0$ for $i\in [d]$. Assume $\de=(\de_1\ldots,\de_d)$ where $\de_i\in \{1,-1\}$ for all $i\in [d]$. Let $H_{\de}$ be the hyperplane whose equation is $\sum_1^d \de_ix_i=1$. The hyperplanes $H_{\de}$ contain the facets of the standard octahedron in $\rr$. Now replace each $H_i$  by  a hyperplane $h_i$ very close to it and, further, each $H_{\de}$ by a hyperplane $h_{\de}$ very close to it. One can of course choose the system
\[
\F=\{h_1,\ldots,h_d\}\cup \{h_{\de}: \mbox{ all } 2^d \pm1 \mbox{ vectors }\de\}
\]
so the intersection of any $d$ hyperplanes from $\F$ is a single point and all of these ${2^d+d \choose d}$ intersection points are distinct. We claim that $\F$ is not in convex position. Assume it is and let $P$ be the polytope with $h\cap P$ a facet of $P$ for every $h\in \F$. Then $P$ must lie in a connected component, say $C$, of the complement of $\bigcup_{h \in \F} h$. The complement of the union of hyperplanes $h_1,\ldots,h_d$ consists of $2^d$ cones, and $C$ is contained in one of them. Each such cone is the intersection of halfspaces of the form $h_i^{\de_i}$ with $\de_i=\pm 1$ where $h_i^+$ and $h_i^-$ are the two halfspaces determined by $h_i$. Assume without loss of generality that $C \subset \bigcap _1^d h_i^+$. But $C$ is disjoint from the hyperplane $h_{\de}$ with $\de=(-1,\ldots,-1)$. So $P\cap h_{\de}\subset C\cap h_{\de}=\emptyset$. So $P\cap h_{\de}$ is not a facet of $P$.\qed

\bigskip
{\bf Acknowledgements.}  We would like to thank two anonymous referees for their helpful comments. 
Research of IB was partially supported by Hungarian National Research grants (no. 131529, 131696, and 133819), and research of GK by the Israel Science Foundation (grant no. 1612/17).

\bigskip


\bigskip

\noindent
Imre B\'ar\'any \\
R\'enyi Institute of Mathematics,\\
13-15 Re\'altanoda Street, Budapest, 1053 Hungary\\
{\tt barany.imre@renyi.hu} and\\
Department of Mathematics\\
University College London\\
Gower Street, London, WC1E 6BT, UK

\medskip
\noindent
Gil Kalai\\
Einstein Institute of Mathematics\\
Hebrew University,
Jerusalem 91904, Israel,\\
{\tt kalai@math.huji.ac.il} and\\
Efi Arazy School of Computer Science,
IDC, Herzliya, Israel

\medskip
\noindent
Attila P\'or\\
Department of Mathematics\\
Western Kentucky University\\
1906 College Heights Blvd. \#11078\\
Bowling Green, KY 42101, USA\\
{\tt  attila.por@wku.edu}\\

\end{document}